# THE AUTOMORPHISM GROUPS
# OF COMPLEX HOMOGENEOUS SPACES

Edward G. Dunne and Roger Zierau

Department of Mathematics
Oklahoma State University
Stillwater, OK 74078

ABSTRACT. If $G$ is a (connected) complex Lie Group and $Z$ is a generalized flag manifold for $G$, then the open orbits $D$ of a (connected) real form $G_0$ of $G$ form an interesting class of complex homogeneous spaces, which play an important role in the representation theory of $G_0$. We find that the group of automorphisms, i.e., the holomorphic diffeomorphisms, is a finite-dimensional Lie group, except for a small number of open orbits, where it is infinite dimensional. In the finite-dimensional case, we determine its structure. Our results have some consequences in representation theory.

§**1.** We determine the automorphism groups for a certain interesting class of complex homogeneous spaces. Denote by $Z$ a generalized flag manifold for a connected complex semisimple Lie group $G$. A real form $G_0$ (which we assume to be connected) of $G$ acts on $Z$ with a finite number of orbits, thus there are always open orbits (cf. [22]). These open orbits play a key role in the representation theory of $G_0$. An open $G_0$-orbit $D$ in $Z$ has a $G_0$-invariant complex structure. The identity component of the group of holomorphic diffeomorphisms of $D$ will be denoted by $Hol(D)$. In the main theorem below we determine $Hol(D)$ for each measurable open orbit (see Definition 2.1). In the case where $D$ is measurable, $D$ carries a $G_0$-invariant (usually) indefinite hermitian metric and we determine its group of hermitian isometries. Generally the open orbits $D$ are non-compact, however, our results include the cases where $G_0$ is a compact real form and $D$ is compact, so $D = Z$. The compact case is contained in [12], [21], [3] and [2], from various points of view.

In general, for a complex manifold $X$, $Hol(X)$ is a (finite-dimensional) Lie group if $X$ is compact and may or may not be a Lie group if $X$ is non-compact. For example, $Hol(\mathbf{C}^n)$ is infinite dimensional. Our main interest is when $G_0$ (so $D$) is non-compact. We give a precise condition for $Hol(D)$ to be a Lie group.

Second author partially supported by N.S.F. Grant DMS 93 03224.

Typeset by $\mathcal{A}_{\mathcal{M}}\mathcal{S}$-TeX

Our main results are contained in the following theorem.

**Main Theorem.** *Suppose $G$ is a connected simple complex Lie group and $D$ is an open measurable $G_0$-orbit in a generalized complex flag manifold for $G$.*

(1) *If there is a $G_0$-equivariant holomorphic fibration of $D$ over the hermitian symmetric space for $G_0$ (and $D$ is not equal to the hermitian symmetric space itself), then $Hol(D)$ is not a finite-dimensional Lie group.*

(2) *If no such fibration exists, then $Hol(D)$ is a Lie group and, except for the cases listed in Table 1.1 below, we have:*

   (a) $Hol(D) = G_0$ *if $G_0$ is non-compact and,*
   
   (b) $Hol(D) = G$ *if $G_0$ is compact.*

| $G_0$ | $Z$ | $Hol(D)$ |
|---|---|---|
| $SO_e(2p, 2q+1)$, $p \neq 0$ | {pure spinors in $\mathbf{C}^{2p+2q+1}$} | $SO_e(2p, 2q+2)/Z_2$ |
| $SO(2n+1)$ | " | $SO(2n+2, \mathbf{C})/Z_2$ |
| $Sp(n, \mathbf{R})$ | $\mathbf{CP}^{2n-1}$ | $SU(n,n)/Z_{2n}$ |
| $Sp(p,q)$, $pq \neq 0$ | " | $SU(2p, 2q)/Z_{2p+2q}$ |
| $Sp(n)$ | " | $SL(2n, \mathbf{C})/Z_{2n}$ |
| $G_{2,\text{split}}$ | quadric in $\mathbf{CP}^6$ | $SO_e(3,4)$ |
| $G_{2,\text{compact}}$ | " | $SO(7, \mathbf{C})$ |

*$D$ is any open orbit in $Z$.*

Table 1.1

(3) *The group of hermitian isometries (i.e., the group of holomorphic diffeomorphisms preserving the hermitian metric) is*

   (a) $Hol(D)$ *if $G_0$ is non-compact and,*
   
   (b) *a compact real form of $Hol(D)$ if $G_0$ is compact.*

In Proposition 3.11 we will see how the case of a semisimple group reduces to the case of simple groups.

The method of proof is to study the Lie algebra of global holomorphic vector fields on $D$ using some standard techniques from representation theory. In most cases, this Lie algebra is just $\mathfrak{g}$. However, in other cases, it is a bigger finite-dimensional Lie algebra $\mathfrak{g}^1$. In each of these cases we find a group $G_0^1$ which has (complexified) Lie algebra $\mathfrak{g}^1$ and has an effective action on $D$.

Our results have several consequences for the representations associated to the open orbits. In the cases listed in Table 1.1 we view $G_0^1$ as acting on $D$. The irreducible representations occurring in Dolbeault cohomology spaces on $D$ extend to (irreducible) representations of $G_0^1$. Also, the results have implications for a space of maximal compact subvarieties of $D$, which in turn plays a role in certain realizations of these representations. This will be discussed in the final section.

We are grateful to T.N. Bailey for helpful discussions regarding Proposition 2.3. and especially to D.A. Vogan for suggesting the approach to Propositions 2.4 and 2.5.

§**2.** Some detailed information is obtained on the Lie algebra of global holomorphic vector fields on $D$. It will follow that $Hol(D)$ is usually finite dimensional and its structure will be narrowed down to a few possibilities. We start with some notation.

Let $G$ be a connected simple complex Lie group. As mentioned in the introduction, the semisimple case can be reduced to the simple case. Fix a generalized complex flag manifold $Z$ for $G$. Then, $Z$ is (biholomorphic to) $G/Q$ where $Q$ is a parabolic subgroup of $G$. We follow the common practice of denoting the Lie algebra of a Lie group by the corresponding gothic letter. Thus, the Lie algebras of $G$ and $Q$ will be denoted by $\mathfrak{g}$ and $\mathfrak{q}$, respectively. A connected real form of $G$ will be denoted by $G_0$, with Lie algebra $\mathfrak{g}_0$. By Theorem 2.6 of [22], $G_0$ acts on $Z$ with a finite number of orbits. Hence we know that open orbits always exist. Fix a Cartan involution $\theta$ and let $K$ (respectively, $K_0$) denote the fixed-point group of $\theta$ in $G$ (respectively, $G_0$). Fix an arbitrary open orbit $D$ in $Z$. For $z_0 \in D$, $\mathrm{Stab}_G(z_0)$ is a parabolic subgroup $Q$. The Levi decompositions are given by $Q = L\overline{U}$ and $\mathfrak{q} = \mathfrak{l} + \overline{\mathfrak{u}}$. By Theorem 4.5 of [22] we may choose $z_0$ so that $\mathfrak{q}$ contains a Cartan subalgebra $\mathfrak{h}$ of $\mathfrak{g}$, so that $\mathfrak{t} = \mathfrak{h} \cap \mathfrak{k}$ is a Cartan subalgebra of $\mathfrak{k}$.

**Definition 2.1.** *The orbit $D$ is* measurable *if and only if $D$ carries a $G_0$-invariant measure.*

By Theorem 6.3 of [22], $D$ is measurable if and only if $\mathfrak{q} \cap \overline{\mathfrak{q}} = \mathfrak{l}$. This is equivalent to the statement that $\mathrm{Stab}_{G_0}(z_0) = L_0$, a real form of $L$. We describe this in terms of the Lie algebras as follows. Let $\Delta = \Delta(\mathfrak{h}, \mathfrak{g})$ be the roots of $\mathfrak{h}$ in $\mathfrak{g}$. Then, there exists $\lambda_0 \in i\mathfrak{t}_0^*$ such that $\Delta(\mathfrak{h}, \mathfrak{l}) = \{\alpha \in \Delta \,|\, \langle \lambda_0, \alpha \rangle = 0\}$ and $\Delta(\mathfrak{h}, \overline{\mathfrak{u}}) = \{\alpha \in \Delta \,|\, \langle \lambda_0, \alpha \rangle < 0\}$. Also, $\overline{\mathfrak{u}}$ is the complex conjugate (with respect to the real form $\mathfrak{g}_0$) of a subalgebra $\mathfrak{u}$ and $\overline{\mathfrak{q}} = \mathfrak{l} + \mathfrak{u}$ is the parabolic opposite to $\mathfrak{q}$. We assume from now on that $D$ is a measurable open $G_0$-orbit in $Z$.

*Remark.* It is the *measurable* open orbits that play an important role in representation theory. If $G_0$ contains a compact Cartan subgroup or if $Z$ is the full flag manifold (so $Q$ is a Borel subgroup), then all open orbits are measurable. Also, if one open $G_0$-orbit in $Z$ is measurable, then all open orbits are measurable. A simple example of a non-measurable open orbit is the (unique) open orbit of $\mathrm{SL}(3, \mathbf{R})$ on $\mathbf{CP}^2$.

The hermitian symmetric spaces $G_0/K_0$ are examples of measurable open orbits.

In this case, writing the Cartan decomposition as $\mathfrak{g} = \mathfrak{k} \oplus \mathfrak{p}$, $\mathfrak{p}$ splits into $\mathfrak{p}_+ \oplus \mathfrak{p}_-$ as representations of $K$. One sees that $G_0/K_0$ is an open orbit in $G/KP_-$. It is well known that $G_0/K_0$ is biholomorphic to a bounded domain in some $\mathbf{C}^N$. Thus the Lie algebra of global holomorphic vector fields is infinite dimensional, so our arguments in this case are a little different from the general case. We will treat this case first.

**Proposition 2.2.** *If $B = G_0/K_0$ is a hermitian symmetric space, then $Hol(B)$ is finite-dimensional.*

*Proof.* This is standard. It is contained in Chapter VIII of [7].

We now give a condition for $Hol(D)$ to be infinite dimensional.

**Proposition 2.3.** *If $B = G_0/K_0$ is of hermitian symmetric type and $D$ is an open orbit with $\mathfrak{p}_- \subsetneq \overline{\mathfrak{u}}$ then the natural fibration $\pi : D \to B$ is holomorphically trivial, in particular $D \cong B \times K_0 \cdot z_0$. Furthermore, the group of holomorphic diffeomorphisms of $D$ is not a finite-dimensional Lie group.*

*Proof.* Recall the following decomposition of $G_0$. In the complex group $G$, $P_+ K P_-$ is a dense open set and $G_0 \subset P_+ K P_-$. Furthermore, the decomposition of $g \in G_0$ as $g = p_+ k p_-$ is unique. Then the Harish-Chandra embedding maps $gK_0 \in B$ to $\xi \in \mathfrak{p}_+$, where $p_+ = \exp(\xi)$. Denote the image of the Harish-Chandra embedding by $\mathcal{B}$. Then $\mathcal{B}$ is a bounded domain biholomorphic to $B = G_0/K_0$.

Now suppose that $D$ has $\mathfrak{p}_- \subset \overline{\mathfrak{u}}$. Then $L_0 \subset K_0$ and $L\overline{U} \subset KP_-$. It is clear that $D = G_0 \cdot z_0 \subset P_+ K \cdot z_0$. The latter is biholomorphic to $\mathfrak{p}_+ \times K_0/(K_0 \cap L_0)$. To see this, note that the map $\Phi : \mathfrak{p}_+ \times K_0/(K_0 \cap L_0) \to P_+ K \cdot z_0$ given by $\Phi(\xi_+, k) = \exp(\xi_+) k \cdot z_0$ is complex analytic. It is clearly onto. It is one-to-one: if $\exp(\xi) k \cdot z_0 = \exp(\xi') k' \cdot z_0$, then by the uniqueness of the Harish-Chandra decomposition, $\exp(\xi) = \exp(\xi')$; so $\xi = \xi'$. Then $k \cdot z_0 = k' \cdot z_0$. We must determine the inverse image of $D$ under $\Phi$. Note that elements of $G_0$ can be written as $g = \exp(\xi) k p_-$ with $\xi \in \mathcal{B}$. Thus, $g \cdot z_0 = \exp(\xi) k \cdot z_0$. As $\Phi$ is bijective, $\Phi^{-1}(D) = \mathcal{B} \times K \cdot z_0$. But $K \cdot z_0 = K_0 \cdot z_0 \cong K_0/(K_0 \cap L_0)$. This proves the first part of the proposition.

As $\mathcal{B}$ is a bounded domain in some $\mathbf{C}^N$, the space of holomorphic functions $f : \mathcal{B} \to \mathfrak{k}$ is infinite dimensional. Now, $K \cdot z_0 = K_0 \cdot z_0$, so $K$ is contained in the automorphism group of $K_0 \cdot z_0$. Thus, for each holomorphic $f : \mathcal{B} \to \mathfrak{k}$, there is an automorphism of $\mathcal{B} \times K_0 \cdot z_0$ defined by $\varphi(\xi, z) = (\xi, (\exp(f(\xi)))z)$. Note that $\mathfrak{k} \neq 0$ since $\mathfrak{p}_- \neq \overline{\mathfrak{u}}$. This provides an infinite-dimensional family in $Hol(D)$. □

Proposition 2.3 is (1) of the main theorem, since there exists a holomorphic fibration $D \to B$ as in (1) if and only if $\overline{\mathfrak{u}}$ properly contains $\mathfrak{p}_+$ or $\mathfrak{p}_-$. Also note that $B$ is simply connected, whence $D \to B$ cannot have discrete fibers.

Assume from now on that $\overline{\mathfrak{u}}$ does not contain $\mathfrak{p}_-$ or $\mathfrak{p}_+$. In these cases, as we will see, $Hol(D)$ is a (finite-dimensional) Lie group. Now consider the Lie algebra of holomorphic vector fields on $D$, that is, the global sections $H^0(D, \mathcal{T})$ of the vector bundle $\mathcal{T}$ of holomorphic vector fields on $D$. Note that $\mathcal{T}$ is the holomorphic homogeneous vector bundle for the representation $\mathfrak{g}/\mathfrak{q}$ of $Q$. Since $G_0$ acts on $D$ and on $\mathcal{T}$, it is clear that $H^0(D, \mathcal{T})$ is a representation of $G_0$. There is considerable machinery available to study representations of this type. See, for example, [17] and [25]. We will very briefly describe the ingredients that we will need. For any finite-dimensional representation $F$ of $Q$, let $\mathcal{F} \to Z$ be the corresponding holomorphic, homogeneous vector bundle on $Z$. By restriction, we obtain a homogeneous vector bundle on $D$ denoted by $\mathcal{F} \to D$. By [25], the Dolbeault cohomology spaces $H^p(D, \mathcal{F})$ are continuous, admissible representations of $G_0$. The Harish-Chandra module (i.e., the subspace of $K_0$-finite vectors in $H^p(D, \mathcal{F})$) is a cohomologically parabolically induced $(\mathfrak{g}, K_0)$-module $\mathcal{R}_\mathfrak{q}^p(F)$. The definition of $\mathcal{R}_\mathfrak{q}^p(F)$ is given in Definition 1, page 432 of [25] (and in §6.3 of [16], with slightly different conventions). It follows immediately that $\mathcal{R}_\mathfrak{q}^0(F) \subset \text{Hom}_\mathfrak{q}(\mathcal{U}(\mathfrak{g}), F)_{K_0-\text{finite}}$. Here, $\mathcal{U}(\mathfrak{g})$ is the universal enveloping algebra of $\mathfrak{g}$ and the subscript '$K_0$-finite' indicates the subspace of $K_0$-finite vectors.

**Proposition 2.4.** *If $E$ is a finite-dimensional representation of $G_0$, then*
$$\text{Hom}_{G_0}(E, H^0(D, \mathcal{T})) \neq 0 \implies \text{Hom}_\mathfrak{q}(E, \mathfrak{g}/\mathfrak{q}) \neq 0.$$

*Proof.* If $\phi : E \to H^0(D, \mathcal{T})$ is a nonzero $G_0$-homomorphism, then the image of $\phi$ lies in $H^0(D, \mathcal{T})$. But $H^0(D, \mathcal{T}) \subset \text{Hom}_\mathfrak{q}(\mathcal{U}(\mathfrak{g}), \mathfrak{g}/\mathfrak{q})$. So $\phi$ defines a nonzero element of $\text{Hom}_\mathfrak{g}(E, \text{Hom}_\mathfrak{q}(\mathcal{U}(\mathfrak{g}), \mathfrak{g}/\mathfrak{q})) \cong \text{Hom}_\mathfrak{q}(E, \mathfrak{g}/\mathfrak{q})$.

**Proposition 2.5.** $H^0(D, \mathcal{T})$ *is finite dimensional and $Hol(D)$ is a Lie group unless $G_0/K_0$ is hermitian symmetric and $\mathfrak{p}_+ \subset \overline{\mathfrak{u}}$ or $\mathfrak{p}_- \subset \overline{\mathfrak{u}}$.*

*Proof.* Fix a positive root system $\Delta^+ \subset \Delta$ so that $\Delta(\mathfrak{u}) \subset \Delta^+$. (We may choose $\Delta^+ = \Delta^+\mathfrak{l} \cup \Delta(\mathfrak{u})$ with $\Delta^+(\mathfrak{l})$ an arbitrary positive system for $\Delta(\mathfrak{l})$.) Suppose $F$ is a finite dimensional $Q = L\overline{U}$ representation with trivial $\overline{U}$ action. Then, viewing $F$ as a $\overline{Q} = LU$ representation with trivial $U$ action, by Lemma 5.15 of [17], $\text{Hom}_\mathfrak{q}(\mathcal{U}(\mathfrak{g}), F)_{L_0 \cap K_0\text{-finite}} \cong \mathcal{U}(\mathfrak{g}) \otimes_{\overline{\mathfrak{q}}} F$. Furthermore, $(\mathcal{U}(\mathfrak{g}) \otimes_{\overline{\mathfrak{q}}} F)_{K_0-\text{finite}}$ is a highest-weight Harish-Chandra module (with respect to $\Delta^+$). By a result of Harish-Chandra (cf. [5], Theorem 1 and Corollary 2, page 761), the only case for which there exist infinite-dimensional highest-weight $(\mathfrak{g}, K)$-modules is when $\Delta^+$ contains $\Delta(\mathfrak{p}_+)$ or $\Delta(\mathfrak{p}_-)$.

Now consider the holomorphic tangent bundle $\mathcal{T}$. The action of $\overline{U}$ on $\mathfrak{g}/\mathfrak{q}$ is not in general trivial, so the above does not apply directly. Instead, we form a filtration $\mathfrak{g}/\mathfrak{q} = F_1 \supset F_2 \supset \cdots \supset F_N \supset 0$ so that $F_i/F_{i+1}$ does have a trivial $\overline{\mathfrak{u}}$ action. (For instance, we could take a composition series for $\mathfrak{g}/\mathfrak{q}$, so that each $F_i/F_{i+1}$ is irreducible as $\mathfrak{q}$-module, in which case the action of $\overline{\mathfrak{u}}$ is necessarily trivial.) As $\text{Hom}_{\mathfrak{q}}(\mathcal{U}(\mathfrak{g}), \cdot)_{(L_0 \cap K_0)\text{-finite}}$ is exact, there is a filtration of $\text{Hom}_{\mathfrak{q}}(\mathcal{U}(\mathfrak{g}), \mathfrak{g}/\mathfrak{q})_{(L_0 \cap K_0)\text{-finite}}$ with quotients $\text{Hom}_{\mathfrak{q}}(\mathcal{U}(\mathfrak{g}), (F_i/F_{i+1}))_{(L_0 \cap K_0)\text{-finite}} \cong (\mathcal{U}(\mathfrak{g}) \otimes_{\overline{\mathfrak{q}}} (F_i/F_{i+1}))$. Now we may conclude that, unless $\Delta^+$ contains one of $\Delta(\mathfrak{p}_\pm)$, $\text{Hom}_{\mathfrak{q}}(\mathcal{U}(\mathfrak{g}), \mathfrak{g}/\mathfrak{q})_{K_0\text{-finite}}$ is finite dimensional.

Now note that as we were free to choose $\Delta^+(\mathfrak{l}) \subset \Delta$, it follows that $\Delta(\mathfrak{p}_\pm) \cap \Delta(\mathfrak{l}) = 0$. Thus, $\Delta^+$ contains one of $\Delta(\mathfrak{p}_\pm)$ if and only if $\Delta(\mathfrak{u})$ contains one of $\Delta(\mathfrak{p}_\pm)$, that is, if and only if $\mathfrak{p}_+ \subset \overline{\mathfrak{u}}$ or $\mathfrak{p}_- \subset \overline{\mathfrak{u}}$.

That $Hol(D)$ is a Lie group now follows from Theorem 3.1 of [9].  □

By Proposition 2.5 we may use Proposition 2.4 to restrict the possibilities for $Hol(D)$. Note, it is clear that

$$\text{Hom}_{\mathfrak{q}}(E, \mathfrak{g}/\mathfrak{q}) \tag{2.6}$$

is non-zero for $E = \mathfrak{g}$ and, of course, that $G_0 \subset Hol(D)$.

**Corollary 2.7.** *If (2.6) is zero for all irreducible representations $E \not\cong \mathfrak{g}$, then $Hol(D) = G_0$ or $G$.*

*Proof.* The Lie algebra of $Hol(D)$ lies between $\mathfrak{g}_0$ and $\mathfrak{g}$. As $\mathfrak{g}$ is simple the only possibilities for $Hol(D)$ are $G_0$ and $G$.  □

We look for finite-dimensional representations $E$ other than $E \cong \mathfrak{g}$ such that (2.6) is non-zero. We start off with a fact about the structure of simple Lie algebras.

**Lemma 2.8.** *If $\mathfrak{g}$ is a simple Lie algebra, then the Weyl group is transitive on roots of a given length. Therefore, there are two possibilities:*
  (1) *There is just one root length and the highest root, which we denote by $\gamma_\ell$, is the unique $\Delta^+$-dominant root.*
  (2) *There are two root lengths and there is a unique dominant long root $\gamma_\ell$ and a unique dominant short root $\gamma_s$.*

*Proof.* This is standard. See, for example, [8], page 53.

**Theorem 2.9.** *If $\mathfrak{g}$ is simple and has just one root length and $\overline{\mathfrak{u}}$ does not contain $\mathfrak{p}_+$ or $\mathfrak{p}_-$, then*

$$Hol(D) = \begin{cases} G_0, & \text{if } G_0 \text{ is non-compact} \\ G, & \text{if } G_0 \text{ is compact.} \end{cases}$$

*Proof.* Suppose $E$ is an irreducible finite-dimensional representation of $G_0$. Then, by Corollary 2.7, it is enough to check that (2.6) is zero unless $E \cong \mathfrak{g}$. Since $\Delta^+$ contains $\Delta(\mathfrak{u})$, the highest-weight space in $E$ is cyclic for $\mathfrak{q} = \mathfrak{l} + \overline{\mathfrak{u}}$. So, a nonzero $\mathfrak{q}$ homomorphism $\phi : E \to \mathfrak{g}/\mathfrak{q}$ maps the highest-weight vector to a (nonzero) weight vector in $\mathfrak{g}/\mathfrak{q}$ (of the same weight). This weight must be dominant. As there is just one root length, Lemma 2.8 says that this highest weight of $E$ must be $\gamma_\ell$. Now Corollary 2.7 applies.

In case $G_0$ is non-compact, $D \neq Z$ and $G$ does not act on $D$. Thus, $Hol(D) = G_0$. When $G_0$ is compact, $D = Z$ and $G$ always acts on $D$. $\square$

When there are two root lengths in $\mathfrak{g}$, then by Lemma 2.8 we have the possibility that the Lie algebra of holomorphic vector fields is $\mathfrak{g}$ or $\mathfrak{g} \oplus E_{\gamma_s}$, where $E_{\gamma_s}$ is the irreducible finite-dimensional representation of $\mathfrak{g}$ with highest weight $\gamma_s$. In connection with this case we have the following lemma.

**Lemma 2.10.** *Let $\mathfrak{g}$ be a simple Lie algebra with two root lengths. Then there is a simple Lie algebra $\mathfrak{g}^1$ containing $\mathfrak{g}$ as a subalgebra so that, as $\mathfrak{g}$-modules, $\mathfrak{g}^1 \cong \mathfrak{g} \oplus E_{\gamma_s}$. The following table lists the possibilities:*

| $\mathfrak{g}$ | $\mathfrak{g}^1$ |
|---|---|
| $B_n$ | $D_{n+1}$ |
| $C_n$ | $A_{2n-1}$ |
| $F_4$ | $E_6$ |
| $G_2$ | $B_3$ |

Table 2.11

*Proof.* This is easily checked by direct calculation.

**§3.** We now consider the simple Lie algebras $\mathfrak{g}$ having two different root lengths. We use Lemma 3.1 below to restrict the possibilities for $Z = G/Q$ for which $Hol(D) \neq G_0$ or $G$. We then treat each of the possible flag manifolds $Z$ separately.

The following lemma will help us determine when (2.6) is zero for $E = E_{\gamma_s}$.

**Lemma 3.1.** *Suppose there are roots $\beta_1, \ldots, \beta_m \in \Delta^+$, not necessarily distinct, so that $\gamma_s - \sum_{j=1}^k \beta_j$ is a root for each $k = 1, \cdots, m$ and $\gamma_s - \sum_{j=1}^m \beta_j$ is a long root in $\mathfrak{u}$. Then $Hom_\mathfrak{q}(E_{\gamma_s}, \mathfrak{g}/\mathfrak{q}) = 0$.*

*Proof.* Suppose $\phi \in Hom_\mathfrak{q}(E_{\gamma_s}, \mathfrak{g}/\mathfrak{q})$ is nonzero. If $v_+$ is a highest weight vector in $E_{\gamma_s}$ then $\phi(v_+)$ is a root vector $X_{\gamma_s}$ of weight $\gamma_s$. Since weights of $E_{\gamma_s}$ all have norm less than or equal to $\| \gamma_s \|$, the long root $\gamma_s - \sum_{j=1}^m \beta_j$ is not a weight. Thus $X_{-\beta_1} \cdots X_{-\beta_m} \cdot v_+ = 0$. On the other hand, $\operatorname{ad}(X_{-\beta_1}) \cdots \operatorname{ad}(X_{-\beta_m}) \cdot X_{\gamma_s} \neq 0$ in $\mathfrak{g}$ since $\gamma_s - \sum_{j=1}^k \beta_j$ is a root for each $k = 1, \cdots, m$. Since the root $\gamma_s - \sum_{j=1}^m \beta_j$ is in $\mathfrak{u}$, $\operatorname{ad}(X_{-\beta_1}) \cdots \operatorname{ad}(X_{-\beta_m}) \cdot X_{\gamma_s} \neq 0$ in $\mathfrak{g}/\mathfrak{q}$. But this is a contradiction; $0 = \phi(X_{-\beta_1} \cdots X_{-\beta_m} \cdot v_+) = \operatorname{ad}(X_{-\beta_1}) \cdots \operatorname{ad}(X_{-\beta_m}) \cdot X_{\gamma_s} \neq 0$. $\square$

For the simple Lie algebras with two root lengths, we number the simple roots as follows:

The fundamental weight corresponding to $\alpha_j$ is denoted by $\lambda_j$. Recall that a weight $\lambda$ determines a parabolic subgroup $Q$ with Lie algebra $\mathfrak{q} = \mathfrak{q}(\lambda)$ by $\Delta(\mathfrak{h}, \mathfrak{q}) = \{\alpha \in \Delta \mid \langle \lambda, \alpha \rangle \leq 0\}$. Each parabolic subalgebra is conjugate to some $\mathfrak{q}(\lambda)$, in fact, is conjugate to one with $\lambda = \sum_{j \in \Phi} \lambda_j$ for some $\Phi \subset \{1, 2, \ldots, \text{rank } \mathfrak{g}\}$.

**Theorem 3.2.** *Hol$(D)$ is $G$ if $G_0$ is compact and is $G_0$ if $G_0$ is non-compact, except possibly in the following cases:*

> *Type $B_n$: $Q$ is determined by $\lambda_n$,*
> *Type $C_n$: $Q$ is determined by $\lambda_1$*
> *Type $G_2$: $Q$ is determined by $\lambda_1$.*

*Proof.* To apply Lemma 3.1 we find, for each simple Lie algebra having more than one root length, roots $\beta_j$ as in the lemma so that $\langle \lambda_j, \gamma_s - \sum_{j=1}^{m} \beta_j \rangle > 0$ for all $\lambda_j$ except those listed in the theorem.

For type $B_n$, $\gamma_s = \sum_{j=1}^{n} \alpha_j$ and $\gamma_s - \alpha_n$ is a long root. (So $m = 1$, $\beta_1 = \alpha_n$ in Lemma 3.1.) Clearly $\langle \lambda_j, \gamma_s - \alpha_n \rangle > 0$ except for $j = n$.

For type $C_n$, $\gamma_s = \alpha_1 + 2\alpha_2 + 2\alpha_3 + \cdots + 2\alpha_{n-1} + \alpha_n$ and $\gamma_s - \alpha_1$ is a long root. (So $m = 1$, $\beta_1 = \alpha_1$ in Lemma 3.1.) Clearly $\langle \lambda_j, \gamma_s - \alpha_1 \rangle > 0$ except for $j = 1$.

For type $G_2$, $\gamma_s = 2\alpha_1 + \alpha_2$ and $\gamma_s - 2\alpha_1(= \alpha_2)$ is a long root. (So $m = 2$, $\beta_1 = \beta_2 = \alpha_2$ in Lemma 3.1.) Clearly $\langle \lambda_2, \gamma_s - 2\alpha_1 \rangle > 0$.

For type $F_4$, $\gamma_s = 2\alpha_1 + 3\alpha_2 + 2\alpha_3 + \alpha_4$ and $\gamma_s - \alpha_2$ is a long root. (So $m = 1$, $\beta_1 = \alpha_2$ in Lemma 3.1.) Clearly $\langle \lambda_j, \gamma_s - \alpha_2 \rangle > 0$ for all $j$. Thus, every parabolic subalgebra of $F_4$ satisfies the condition of Lemma 3.1. □

We now prove our main theorem by looking at each flag manifold listed in Theorem 3.2.

**Type $B_n$.** Let $p, q$ be positive integers such that $p + q = n$. The complex group $G = \mathrm{SO}(2n+1, \mathbf{C})$ is defined by the symmetric form on $\mathbf{C}^{2n+1}$:

$$(w, z) = \sum_{j=1}^{2p} w_j z_j - \sum_{j=2p+1}^{2p+2q+1} w_j z_j.$$

The relevant real form $G_0$ of $B_n$ is the connected component of the isometry group of $(\,,\,)$ restricted to $\mathbf{R}^{2n+1}$. That is, $G_0 \cong \mathrm{SO}_e(2p, 2q+1)$. The complex flag manifold $Z$ is the space of all maximal isotropic subspaces of $\mathbf{C}^{2n+1}$, also known as the space of pure spinors. The set of $(p+q)$-planes $\zeta \in Z$ such that the hermitian form $\langle w, z\rangle := (w, \bar z)$ restricted to $\zeta$ has signature $(p, q)$ contains two open $G_0$-orbits. They are $D_\pm = G_0.z_0^\pm$, where $z_0^\pm = \mathrm{span}_{\mathbf{C}}\{e_1 \pm ie_2, e_3 + ie_4, \ldots, e_{2n-1} + ie_{2n}\}$. These are the only open orbits of $G_0$ in $Z$.

Let $G^1$ be the isometry group of $(w, z)^1 = \sum_{j=1}^{2p} w_j z_j - \sum_{j=2p+1}^{2p+2q+2} w_j z_j$ on $\mathbf{C}^{2n+2}$. Let $G_0^1$ be the connected component of the isometry group of $(\,,\,)^1$ restricted to $\mathbf{R}^{2n+2}$. Then, $G_0^1 \cong \mathrm{SO}_e(2p, 2q+2)$. The space of maximal isotropic subspaces in $\mathbf{C}^{2n+2}$ has two connected components (since the dimension is even). Take $Z^1$ to be the component containing $z_0^1 = \mathrm{span}_{\mathbf{C}}\{e_1 + ie_2, \ldots, e_{2n+1} + ie_{2n+2}\}$. (The other component contains $\mathrm{span}_{\mathbf{C}}\{e_1 + ie_2, \ldots, e_{2n-1} + ie_{2n}, e_{2n+1} - ie_{2n+2}\}$.) Then $Z^1$ is a flag manifold for $G^1$. Again, there are two open orbits $D_\pm^1$ consisting of planes of signature $(p, q+1)$ with respect to the hermitian form $\langle w, z\rangle^1 := (w, \bar z)^1$.

Define a map $\pi : Z^1 \to Z$ by $\pi(\zeta) = \zeta \cap (\mathbf{C}^{2n+1} \times \{0\})$. Note that $\pi(\zeta) \in Z$ as it is clearly isotropic and $\dim(\pi(\zeta)) = n$. It is also clear that $\pi$ is $G_0$-equivariant. By equivariance, $\pi$ is onto. To see that $\pi$ is one-to-one it is enough to check that $\pi^{-1}(z_0) = \{z_0^1\}$. Suppose $\pi(\zeta) = \zeta \cap (\mathbf{C}^{2n+1} \times \{0\}) = z_0$. Then $z_0 \subset \zeta$ and there is $v \in \zeta$ so that $v \in z_0^{\perp_{(\,,\,)}} \cap z_0^{\perp_{\langle\,,\,\rangle}} = (z_0 + \overline{z_0})^\perp$ and $\zeta = z_0 \oplus \mathbf{C} \cdot v$. But $(z_0 + \overline{z_0})^\perp = \mathrm{span}_{\mathbf{C}}\{e_{2n+1}, e_{2n+2}\}$. As $v$ is isotropic it is either $e_{2n+1} + ie_{2n+2}$ or $e_{2n+1} - ie_{2n+2}$. Only $v = e_{2n+1} + ie_{2n+2}$ gives $\zeta \in Z^1$, so $\zeta = z_0^1$.

We may conclude that $G$ acts transitively on $Z^1 \cong Z$ and $G_0$ acts transitively on $D_\pm^1 \cong D_\pm$. Equivalently, $G^1$ and $G_0^1$ act on $Z$ and $D_\pm$, respectively. It also follows that the compact real form acts on the $G_0$-orbit $Z$.

**Type $C_n$.** Let $\omega(w, z) = \sum_{j=1}^{n}(w_j z_{n+j} - w_{n+j} z_j)$ be the standard symplectic form on $\mathbf{C}^{2n}$. Then $G = \mathrm{Sp}(n, \mathbf{C})$ is the complex group preserving $\omega$. The complex flag manifold under consideration is $Z = \{\omega\text{-isotropic lines in } \mathbf{C}^{2n}\}$, which is just $\mathbf{CP}^{2n-1}$, since any *line* is automatically isotropic. There are two families of real forms: $\mathrm{Sp}(n, \mathbf{R})$ and $\mathrm{Sp}(p, q)$.

Define $G_0 = \mathrm{U}(n, n) \cap G$, where $\mathrm{U}(n, n)$ is the isometry group of $\langle w, z\rangle = \sum_{j=1}^{n}(w_j \overline{z_j} - w_{n+j}\overline{z_{n+j}})$. Then $G_0 \cong \mathrm{Sp}(n, \mathbf{R})$. There are two open $G_0$-orbits, $D_\pm$: the positive lines and the negative lines. But, it is clear from Witt's theorem that $G_0^1 = \mathrm{SU}(n, n)$ acts transitively on $D_\pm$.

For the other real forms, take $G_0 = \mathrm{U}(2p, 2q) \cap G$, where $\mathrm{U}(2p, 2q)$ is the isometry group of $\langle w, z \rangle = \sum_{j=1}^{p}(w_j \overline{z_j} + w_{n+j}\overline{z_{n+j}}) - \sum_{j=p+1}^{p+q}(w_j \overline{z_j} + w_{n+j}\overline{z_{n+j}})$, where $p + q = n$. Then $G_0 \cong \mathrm{Sp}(p, q)$, (cf. [7], p. 445). For $pq \neq 0$, $G_0$ has two open orbits consisting of positive and negative lines. Also, $G_0^1 = \mathrm{SU}(2p, 2q)$ acts transitively on each of these orbits. If $pq = 0$, then $G_0$ is compact and acts transitively on $Z = \mathbf{CP}^{2n-1}$. The unitary group $\mathrm{SU}(2n)$ also acts transitively on $Z$.

**Type $\mathrm{G}_2$.** Let $G$ be the complex group of type $\mathrm{G}_2$ and let $G^1$ be the complex group of type $\mathrm{B}_3$. Let $G_0$ and $G_0^1$ be the split real forms of $G$ and $G^1$, respectively. We first recall how $G_0 \subset G_0^1$ and $G \subset G^1$. The split real form $G_0$ of type $\mathrm{G}_2$ is the automorphism group of the split octonions (i.e., the Cayley numbers), $\tilde{\mathbf{O}}$. There is a natural inner product of signature $(4, 4)$ defined on $\tilde{\mathbf{O}}$. It is easy to see that any automorphism must preserve the inner product. Since an automorphism must fix $1 \in \tilde{\mathbf{O}}$, it must preserve $\mathrm{Im}\,\tilde{\mathbf{O}} := (\mathbf{R} \cdot 1)^\perp$, the *pure imaginary octonions*. Moreover, an automorphism is completely determined by its restriction to $\mathrm{Im}\,\tilde{\mathbf{O}}$. We see, then, that $G_0 \subseteq \mathrm{SO}_e(3, 4) \cong G_0^1$. Complexifying, we obtain a complex symmetric form $(\,,\,)_{\mathbf{C}}$ on $(\mathrm{Im}\,\tilde{\mathbf{O}})_{\mathbf{C}}$. Similarly, $G$ sits inside the isometry group $G^1$ of $(\,,\,)_{\mathbf{C}}$. There is a corresponding hermitian form defined by $\langle w, z \rangle := (w, \bar{z})_{\mathbf{C}}$.

Consider the flag manifold for $G^1$ defined by $Z = \{\text{isotropic lines in } (\mathrm{Im}\,\tilde{\mathbf{O}})_{\mathbf{C}} \cong \mathbf{C}^7\}$. This is a flag manifold for $G$ under the action of $G$ as a subgroup of $G^1$. The proof of this is similar to the proof of Claim 3.5 below and is sketched at the end of this section. We will see that $G_0^1 = \mathrm{SO}_e(3, 4)$ acts on $Z$ with two open orbits $D_\pm$ and that $G_0$ acts transitively on both.

The action of $G_0^1 \cong \mathrm{SO}_e(3, 4)$ has two open orbits, $D_\pm$, respectively, the positive and negative lines (with respect to $\langle\,,\,\rangle$) in $Z$. We verify this as follows. Write $z = x + iy \in (\mathrm{Im}\,\tilde{\mathbf{O}})_{\mathbf{C}} = \mathrm{Im}\,\tilde{\mathbf{O}} + i\,\mathrm{Im}\,\tilde{\mathbf{O}}$ with $x, y \in \mathrm{Im}\,\tilde{\mathbf{O}}$. Then $z$ is positive and isotropic if and only if:

$$(3.4) \qquad (x, x) = (y, y) > 0 \text{ and } (x, y) = 0\,.$$

Suppose $z' = x' + iy'$ is another positive isotropic vector. By scaling $z$, we may assume $(x, x) = (x', x')$ and $(y, y) = (y', y')$. Now, Witt's theorem says that there is an isometry of $\mathrm{Im}\,\tilde{\mathbf{O}} \cong \mathbf{R}^{3,4}$ taking $x$ to $x'$ and $y$ to $y'$, i.e., taking $z$ to $z'$. One can check that this isometry can be chosen to lie in $\mathrm{SO}_e(3, 4)$.

**Claim 3.5.** $G_0$ *acts transitively on* $D_\pm$.

We use the following lemma, which follows easily from the development in Chapter 6 of [6].

**Lemma 3.6.** *Suppose $A$, $A' \subset \tilde{\mathbf{O}}$ are normed subalgebras both isomorphic to either the quaternions $\mathbf{H}$ or the $2 \times 2$ real matrices $M_2(\mathbf{R})$, then any isomorphism $A \to A'$ extends to an automorphism of $\tilde{\mathbf{O}}$.*

*Proof.* Choose $\epsilon \in A^\perp$ and $\epsilon' \in (A')^\perp$ with $(\epsilon, \epsilon) = -1 = (\epsilon', \epsilon')$. Then, by Lemma 6.15 of [6], $A^\perp = A\epsilon$ and $\tilde{\mathbf{O}} = A \oplus A\epsilon$. Moreover, multiplication in $\tilde{\mathbf{O}}$ is given by

$$(3.7) \qquad (a + b\epsilon)(c + d\epsilon) = (ac + \bar{d}b) + (da + b\bar{c})\epsilon,$$

where the bar denotes conjugation in $A$: the usual conjugation for $\mathbf{H}$ and

$$\overline{\begin{pmatrix} \alpha & \beta \\ \gamma & \delta \end{pmatrix}} = \begin{pmatrix} \delta & -\beta \\ -\gamma & \alpha \end{pmatrix}$$

for $M_2(\mathbf{R})$. Similarly for $\tilde{\mathbf{O}} = A' \oplus A'\epsilon'$. Now suppose $f : A \to A'$ is an isomorphism. Then, as one can readily check, $F(a + b\epsilon) := f(a) + f(b)\epsilon'$ is an automorphism of $\tilde{\mathbf{O}}$. $\square$

*Remark 3.8..* In case $A = A'$, $F_{\epsilon'}(a + b\epsilon) := f(a) + f(b)\epsilon'$ extends $f$ to an automorphism of $\tilde{\mathbf{O}}$ for any $\epsilon' \in A^\perp$ with $(\epsilon', \epsilon') = -1$. In fact, any automorphism of $\tilde{\mathbf{O}}$ that preserves $A$ must send $\epsilon$ to such an $\epsilon'$. Therefore, all extensions of $f$ are of this form.

*Proof of Claim 3.5.* If $z$ and $z'$ are two elements in $D_+$ with $z = x + iy$ and $z' = x' + iy'$, then both decompositions must satisfy (3.4). By Proposition 6.40 of [6], the subalgebras $A$ and $A'$ of $\tilde{\mathbf{O}}$ generated by $\{x, y\}$ and $\{x', y'\}$ are isomorphic to $\mathbf{H}$, say, $\phi : A \to \mathbf{H}$ and $\phi' : A' \to \mathbf{H}$. Then, $f = \phi^{-1} \circ \phi' : A' \to A$ is an isomorphism from $A'$ to $A$.

Now, both $x$, $y$ and $f(x')$, $f(y')$ satisfy (3.4). By rescaling $z'$, if necessary, we may assume $(x, x) = (f(x'), f(x'))$. Also, $\mathrm{Aut}(\mathbf{H})$ consists solely of the conjugations by unit quaternions, which is the map $\mathrm{SU}(2) \to \mathrm{SO}(3)$, acting on $\mathrm{Im}(\mathbf{H}) \cong \mathbf{R}^3$. By Witt's theorem, there is an orthogonal map $h$ (i.e. an automorphism of $\mathbf{H}$) sending $f(x')$ to $x$ and $f(y')$ to $y$. Then, $h \circ f$ extends to an automorphism $g$ of $\tilde{\mathbf{O}}$, by Lemma 3.6. Thus, we have an element $g \in G_0$ with $g(z') = z$.

For $z, z' \in D_-$, the argument is similar. Now, however, the subalgebras $A$ and $A'$ generated by $\{x, y\}$ and $\{x', y'\}$ are isomorphic to $M_2(\mathbf{R})$. Again, there exists an isomorphism $f : A' \to A$. The group $\mathrm{Aut}(M_2(\mathbf{R}))$ consists solely of conjugations by elements of $\mathrm{SL}(2, \mathbf{R})$, which is the action of $\mathrm{SO}_e(1, 2)$ on $\mathbf{R}^{1,2} \cong M_2(\mathbf{R}) \cap \mathrm{Im}\,\tilde{\mathbf{O}}$. Again, by Witt's theorem, there is an isometry $h$ of $\mathbf{R}^{1,2}$ such that $h(f(x')) = x$ and $h(f(y')) = y$. One can check that $h$ may be taken to lie in $\mathrm{SO}_e(1, 2) \cong \mathrm{SL}(2, \mathbf{R})$. Now, by Lemma 3.6, $h \circ f$ extends to $g \in \mathrm{Aut}(\tilde{\mathbf{O}})$. Thus, we have an element $g \in G_0$ with $g \cdot z' = z$. This completes the proof of Claim 3.5. $\square$

*Remark 3.9.* As homogeneous spaces, the two open orbits are $D_+ \cong G_0/\mathrm{U}(2)$ and $D_- \cong G_0/(\mathrm{SL}(2, \mathbf{R}) \times \mathrm{U}(1))$. This is a consequence of Lemma 3.6 and Remark 3.8, as we now demonstrate.

For $D_+$, write $A = \mathbf{H} = \bigoplus_{j=1}^{4} \mathbf{R}e_j$ with $e_1$ = identity and $e_2, e_3, e_4$ units (i.e. $(e_j, e_j) = 1$) with $e_2 e_3 = e_4$. Suppose $g \in G_0$ fixes $e_2 + ie_3 \in D_+$. Then $g$ fixes both $e_2$ and $e_3$. As $g$ is an automorphism, it must also fix $e_4$, hence all of $A$. By Remark 3.8, the automorphisms fixing $A$ are of the form $F_{\alpha\epsilon}$, where $\alpha$ is an element of $A$ of length one. But, since $A \cong \mathbf{H}$, the length-one elements form a group isomorphic to $\mathrm{SU}(2)$. The stabilizer of the line $\mathbf{C} \cdot (e_2 + ie_3)$ also includes the scalars. Therefore, $\mathrm{Stab}_{G_0}(e_2 + ie_3) \cong \mathrm{U}(2)$.

For the orbit $D_-$, take $A = M_2(\mathbf{R})$. We pick our base point to be the isotropic vector $x + iy$ with $x = \begin{pmatrix} 0 & 1 \\ 1 & 0 \end{pmatrix}$ and $y = \begin{pmatrix} 1 & 0 \\ 0 & -1 \end{pmatrix}$. Again, if $g \in G_0$ fixes $x + iy$, then $g$ fixes $A$. By Remark 3.8, the automorphisms of $\tilde{\mathbf{O}}$ fixing $A$ are all of the form $F_{\beta\epsilon}$ where $\beta \in A$ and $(\beta, \beta) = 1$, i.e., $\det(\beta) = 1$. Thus, the group of automorphisms fixing the vector $x + iy$ is $\mathrm{SL}(2, \mathbf{R})$. As scalars also fix the line $\mathbf{C} \cdot (x + iy)$, we have $\mathrm{Stab}_{G_0}(x + iy) \cong \mathrm{SL}(2, \mathbf{R}) \times \mathrm{U}(1)$.

We conclude the type $\mathrm{G}_2$ case by sketching a proof that $Z$ is a flag manifold for $G$. For this it is enough to show that a compact real form of $\mathrm{G}_2$ acts transitively on $Z$. To do this we use a slightly different realization for $G$. A compact real form of $\mathrm{G}_2$ is the automorphism group of $\mathbf{O}$, where $\mathbf{O}$ is the nonsplit octonians. There is a positive definite symmetric form on $\mathbf{O}$ preserved by the automorphism group, $\mathrm{Aut}(\mathbf{O})$. As above $\mathrm{Aut}(\mathbf{O})$ preserves $\mathrm{Im}\,\mathbf{O} = (\mathbf{R} \cdot 1)^\perp$, thus $\mathrm{Aut}(\mathbf{O}) \subseteq \mathrm{SO}(7)$. Now complexify the form to obtain a symmetric form on $(\mathrm{Im}\,\mathbf{O})_\mathbf{C}$. The isometry group is conjugate to $G^1$ (as the form is equivalent to the complex form arising from the split octonions). The corresponding flag manifold is the space of isotopic lines in $\mathrm{Im}\,\mathbf{O}$.

Now let $z = x + iy$ and $z' = x' + iy'$ be isotropic vectors in $(\mathrm{Im}\,\mathbf{O})_\mathbf{C}$. Then (3.4) holds for both $z$ and $z'$. The subalgebras $A$, $A'$ generated by $\{x, y\}$ and $\{x', y'\}$ are both isomorphic to $\mathbf{H}$. Now argue as in the case $D_+$ using the fact, analogous to Lemma 3.6, that an isomorphism $A \to A'$ extends to an automorphism of $\mathbf{O}$.

Part (3) of the main theorem now follows immediately since the groups $G_0$ or $G_0^1$ preserve the metric and the corresponding complex groups do not. This completes the proof of the main theorem.

Now suppose that $G$ is semisimple. Write $\mathfrak{g} = \mathfrak{g}_1 \oplus \cdots \oplus \mathfrak{g}_N$ with $\mathfrak{g}_1, \ldots, \mathfrak{g}_N$ simple and let $G_i$ be the subgroups of $G$ corresponding to the $\mathfrak{g}_i$. A flag manifold for $G$ is the product of flag manifolds $Z_i \cong G_i/Q_i$ for the $G_i$. A measurable open $G_0$-orbit $D$ is the product of measurable open $G_{i,0}$-orbits $D_i$ in the $Z_i$. Recall that by Proposition 2.3, writing $\mathfrak{q}_i = \mathfrak{l}_i + \overline{\mathfrak{u}}_i$, if

(3.10) $$\mathfrak{p}_{i,-} \subset \overline{\mathfrak{u}}_i \text{ or } \mathfrak{p}_{i,+} \subset \overline{\mathfrak{u}}_i$$

then $Hol(D_i)$ is not a (finite-dimensional) Lie group. Thus, if at least one of the factors $D_i$ satisfies (3.10) then $Hol(D)$ is not a Lie group.

**Proposition 3.11.** *If no factor $D_i$ satisfies (3.10) or is a hermitian symmetric space then $Hol(D) \cong Hol(D_1) \times \cdots \times Hol(D_N)$.*

*Proof.* By the same argument as in 2.5, $Hol(D)$ is finite dimensional. We look for representations $E$ such that $\operatorname{Hom}_{\mathfrak{q}}(E, \mathfrak{g}/\mathfrak{q}) \neq 0$. The only such irreducible representations $E$ are $E \cong \mathfrak{g}_i$ and $E \cong \mathfrak{g}_i^1$, the $\mathfrak{g}_i^1$ occuring when $Z_i$ appears in Table 1.1. It follows that $Hol(D) = Hol(D_1) \times \cdots \times Hol(D_n)$.

§**4.** The main theorem has several consequences in representation theory. As mentioned earlier, interesting representations are those occurring as Dolbeault cohomology of a line bundle over an open orbit $D$. These are the representations associated to elliptic co-adjoint orbits by the orbit method. See for example [14], [15], [16] and [18] for some of the important properties of these representations. Suppose we have a holomorphic homogeneous line bundle $\mathcal{L}_\chi \to G_0/L_0$ corresponding to the character of $L_0$ whose differential is $\chi \in it_0^*$. If

$$(4.1) \qquad \langle \chi + \rho, \beta \rangle < 0, \text{ for all } \beta \in \Delta(\mathfrak{u})$$

(where $\rho$ is half the sum of the positive roots) then the Dolbeault cohomology space $H^p(G_0/L_0, \mathcal{L}_\chi) = 0$ for $p \neq s = \dim_C(K_0 \cdot z_0)$ and $H^s(G_0/L_0, \mathcal{L}_\chi)$ is an irreducible admissible representation. When $G_0$ is non-compact the main theorem gives the cases where $Hol(D) = G_0^1 \not\supseteq G_0$. From the case-by-case descriptions given in Section 3, it is clear that $D \cong G_0/L_0 \cong G_0^1/L_0^1$ is an open orbit in a flag manifold for $G^1$. The character defining $\mathcal{L}_\chi$ extends uniquely to a character of $L_0^1$ whose differential we will call $\chi^1$. Thus the bundles $\mathcal{L}_\chi \to G_0/L_0$ and $\mathcal{L}_{\chi^1} \to G_0^1/L_0^1$ are holomorphically and $G_0$-equivariantly equivalent, implying the Dolbeault cohomology spaces are $G_0$-isomorphic. Thus, in the range where (4.1) holds, the representations $H^s(G_0^1/L_0^1, \mathcal{L}_\chi)$ are irreducible representations of $G_0^1$ which remain irreducible when restricted to $G_0$. Put differently, the irreducible representations $H^s(G_0/L_0, \mathcal{L}_\chi)$ of $G_0$ extend to representations of $G_0^1$. This is a somewhat rare phenomenon. Wolf [20] obtained some general results on the irreducibility of a representation when restricted to a subgroup, which have some overlap with this application of our main theorem. Also, see [10] for results on restricting representations of this type. More recently, A. Dvorsky has a number of new results along these lines.

A well-known and important example of the phenomenon discussed above occurs in quantum electrodynamics. A family of "massless" representations of the de Sitter group (SO(2,3)) extend to representations of the conformal group (SO(2,4)), cf. [1]. As these representations can be realized naturally in Dolbeault cohomology (for the open orbits in the first entry of Table 1.1 for $p = 1$ and $q = 1$), the geometric explanation for this extension is an instance of our main theorem.

When $G_0$ is compact, for arbitrary $\chi$, $H^p(G_0/L_0, \mathcal{L}_\chi) \cong H^p(G_0^1/L_0^1, \mathcal{L}_{\chi^1})$ as $G_0$ representations for all $p$. By the Bott-Borel-Weil theorem we see that certain representations of $G_0^1$ remain irreducible when restricted to $G_0$. The following table gives the highest weights $\chi$ and $\chi^1$ of these representations.

| $\mathfrak{g}^1$ | $\chi^1$ | $\mathfrak{g}$ | $\chi$ |
|---|---|---|---|
| $D_{n+1}$ | $a\lambda_{n+1}$ | $B_n$ | $a\lambda_n$ |
| $A_{2n-1}$ | $a\lambda_1$ | $C_n$ | $a\lambda_1$ |
| $B_3$ | $a\lambda_1$ | $G_2$ | $a\lambda_1$ |

$a \in Z_+$

Table 4.2

Our results also have implications for the study of deformations of maximal compact subvarieties in $D$, as in [4], [19], [24]. In particular, let $V_0$ be the compact subvariety $K_0 \cdot z_0$ in $D = G_0/L_0$, where $G_0$ is non-compact and $D$ is, as usual, an open orbit in a generalized flag manifold for $G$. The cycle space $M_D$ studied in [19] is the space of translations of $V_0$ by elements of $G$ which remain in $D$: $M_D = \{g \cdot V \subset D \mid g \in G\}$. The space $M_D$ plays a key role in constructing a transform, often called a Penrose transform, for the representations $H^s(D, \mathcal{L}_\chi)$. In most cases, $M_D$ has the dimension that one computes using Kodaira-Spencer theory (cf. [11]). However, in some cases, the Kodaira-Spencer theory predicts a larger space of deformations. For some of these cases (but not all), our results explain the discrepancy by demonstrating that there is a larger group of automorphisms $G_0^1$ acting on $D$. Compare with [13].

For example, when $\mathfrak{g}$ is of type $C_n$ and $G_0 = \text{Sp}(n, \mathbf{R})$ we have seen that the space of positive lines in $\mathbf{CP}^{2n-1}$ is an open orbit. Then $V_0 \cong \mathbf{CP}^{n-1}$. Kodaira-Spencer theory predicts that the (infinitesimal) deformations of $V_0$ in $D$ come from $H^0(V_0, \mathcal{N})$, where $\mathcal{N}$ is the holomorphic conormal bundle of $V_0$ in $D$. In this example, $H^0(V_0, \mathcal{N}) \cong \{\text{symmetric } n \times n \text{ matrices}\} \oplus \{\text{skew-symmetric } n \times n \text{ matrices}\}$ as representations of $K_0$. Let's denote the deformation space of [19] for the groups $G_0 = \text{Sp}(n, \mathbf{R})$ and $G_0^1 \cong \text{SU}(n, n)$ acting on $D$ by $M_D$ and $M_D^1$, respectively. One can check that in this case $M_D \cong G_0/K_0$ and $M_D^1 \cong G_0^1/K_0^1$ (cf. [24]). As $\dim_{\mathbf{C}}(G_0^1/K_0^1) = \dim(G_0/K_0) + \frac{n(n-1)}{2}$, the deformation space of Kodaira-Spencer is strictly larger than $M_D$ in dimension. The extra $\frac{n(n-1)}{2}$ dimensions come from the skew matrices $B$ by the action of

$$\exp\begin{pmatrix} 0 & B \\ \overline{B} & 0 \end{pmatrix} \in U(n,n) \smallsetminus \text{Sp}(n, \mathbf{R})$$

on $V_0$. Note that the infinitesimal deformations predicted by Kodaira-Spencer theory all come from $M_D^1$.

As a final remark, note that the list in [2] of nilpotent co-adjoint orbits 'sharing' an orbit with a bigger group has a lot in common with our list. The open orbits we are considering here are $G_0$-equivariantly biholomorphic to the elliptic co-adjoint orbits. It would be interesting to understand the connection between our list and that of [2] in terms of nilpotent orbits as limits of elliptic orbits.